\documentclass[12pt]{article}

\usepackage{latexsym}
\usepackage{amsmath}
\usepackage{amssymb}
\usepackage{stackrel}
\usepackage{color}
\usepackage{graphicx}

\usepackage{hyperref}

\usepackage[left=2cm,right=2cm,
    top=2cm,bottom=2cm,bindingoffset=0cm]{geometry}

\newtheorem{lemma}{Lemma}

\newtheorem{theorem}{Theorem}
\newtheorem{corol}{Corollary}
\newcommand{\C}{\mathbb{C}}

\newcommand{\Z}{\mathbb{Z}}

\begin{document}

\begin{center}
{\large\bf Small divisors in the problem of the convergence of generalized power series solutions of $q$-difference equations}
\end{center}

\begin{center}
{\large R.\,Gontsov, I.\,Goryuchkina, A.\,Lastra}
\end{center}


\section{Introduction}

In this work we consider a $q$-difference equation
\begin{equation}\label{e1}
F(z,y,\sigma y,\sigma^2 y,\ldots,\sigma^{n}y)=0,
\end{equation}
where $F=F(z,y_0,y_1,\ldots,y_n)$ is a polynomial and $\sigma$ stands for the dilatation operator
$$
\sigma:\;y(z)\mapsto y(qz),
$$
$q\ne0,1$ being a fixed complex number. We study the question of the convergence of its {\it generalized} formal power series solutions
$y=\varphi$ of the form
\begin{equation}\label{eq1}
  \varphi=\sum_{j=0}^{\infty}c_jz^{\lambda_j},\qquad c_j,\lambda_j\in\mathbb{C},
\end{equation}
where $c_0\neq 0$ and the sequence of the exponents $\lambda_j$ possesses the following two properties:
\begin{itemize}
\item[(i)] ${\rm Re}\,\lambda_j\leqslant\hbox{Re}\,\lambda_{j+1}$ for all $j\geqslant0$,
\item[(ii)] $\lim_{j\to\infty}{\rm Re}\,\lambda_j=+\infty$.
\end{itemize}

We note that the conditions (i), (ii) make the set of all generalized formal power series an algebra over $\C$. The definition of the dilatation operator extends naturally to this algebra after fixing the value of $\ln q$ by the condition $0\leqslant{\rm arg}\,q<2\pi$:
$$
\sigma\Bigl(\sum_{j=0}^{\infty}c_jz^{\lambda_j}\Bigr)=\sum_{j=0}^{\infty}c_jq^{\lambda_j}z^{\lambda_j}.
$$
Thus the notion of a generalized formal power series solution of (\ref{e1}) is correctly defined in view of the above remarks: such a series $\varphi$
is said to be a {\it formal solution} of (\ref{e1}) if the substitution of $\varphi$ into the polynomial $F$ leads to a generalized power series with zero coefficients.

Formal solutions \eqref{eq1} generalize classical power series solutions of the form $\sum_{j=0}^{\infty}c_jz^j$ whose
convergence is widely studied yet: there are two principally different cases, that of $|q|\neq 1$ and of $|q|=1$. In the case of $|q|\neq 1$ the sufficient condition of convergence of a formal solution $\sum_{j=0}^{\infty}c_jz^j$ of \eqref{e1}, see \cite{zh}, \cite{lz}, is similar to the corresponding Malgrange condition \cite{m} for the differential case, while the study of convergence in the case of $|q|=1$ is more subtle since the small divisors phenomenon arises in this situation. Below we recall some of the most known results concerning the problem of the convergence of classical power series solutions of \eqref{e1} in the case of $|q|=1$, $q$ not being a root of unity.

We start with Siegel's study of the linearization of a diffeomorphism $f(z)=qz+\dots$ of $(\mathbb{C},0)$. One says that $f$ is linearizable if it belongs to the conjugacy class of the rotation $\tilde f(z)=qz$ in the group $G$ of diffeomorphisms of $(\mathbb{C},0)$. The existence of a linearizing map is equivalent to the existence and convergence of a formal power series solution of the Schr\"{o}der equation (a first order $q$-difference equation)
$$
y(qz)=f(y(z)),\qquad y(0)=0.
$$
Such a solution does exist and converge if
\begin{equation}\label{siegel}
|q^j-1|>c\,j^{-\nu}\qquad \mbox{ for all }j\in\mathbb N,
\end{equation}
where $c$ and $\nu$ are some positive constants \cite{Siegel}. This means that under the assumption \eqref{siegel} the set $G_q\subset G$ of diffeomorphisms whose derivative at $0$ equals $q$, is a conjugacy class in $G$, that is, any $f\in G_q$ is linearizable. Representing $q$ in the form $q=e^{2\pi{\rm i}\omega}$, $\omega$ being irrational, one obtains that \eqref{siegel} is equivalent to
\begin{equation}\label{siegeladd}
|j\omega-m|>c\,j^{-\nu}\qquad \mbox{ for all }j\in{\mathbb N}, m\in\mathbb Z.
\end{equation}
Considering the sequence of the convergents $p_k/q_k$ of the continued fraction for $\omega$, one also establishes the equivalence of \eqref{siegeladd} to the condition
\begin{equation}\label{siegelconv}
\ln q_{k+1}=O(\ln q_k)\qquad \mbox{ for }k\to\infty
\end{equation}
(see, for example, Theorems 9, 13, 16 in \cite{Khinchin}). Further the condition \eqref{siegelconv} of Siegel, sufficient for $G_q$ to be a conjugacy class in $G$, was weakened by R\"ussmann \cite{Rus} and Bruno \cite[Ch. II, Th. 6]{Bruno} who have established a less restrictive sufficient condition
$$
\sum_{k=0}^\infty\frac{\ln q_{k+1}}{q_k}<+\infty.
$$
Finally Yoccoz \cite{Yoccoz} has proved that the Bruno--R\"ussmann condition on $q$ is not only sufficient but necessary as well for that every element of $G_q$ to be linearizable (for $G_q$ to be a conjugacy class in $G$).

For the  general  equation \eqref{e1} of an arbitrary order $n$, the first results for $q=e^{2\pi{\rm i}\omega}$, $\omega$ being irrational, were obtained by B\'ezivin and concern the {\it linear} case, that of (\ref{e1}) with
$$
F=a_0(z)y_0+a_1(z)y_1+\ldots+a_n(z)y_n+b(z),
$$
where the $a_k$'s and $b$ are analytic at the origin. Presenting each $a_k(z)\equiv F'_{y_k}$ in the form
\begin{equation}\label{lin}
a_k(z)=A_kz^d+B_kz^{d+1}+\ldots,
\end{equation}
$d\in{\mathbb Z}_+$ being the same for all $k=0,1,\ldots,n$ and at least one of the $A_k$'s being non-zero, one defines a non-zero polynomial
$Q(\xi)=A_n\xi^n+\ldots+A_1\xi+A_0$. If for each root $a$ of the polynomial $(\xi-1)Q(\xi)$ there holds the estimate
\begin{equation}\label{e2}
|q^j-a|> cj^{-\nu}\qquad \mbox{ for all }j\in\mathbb N,
\end{equation}
or, equivalently,
\begin{equation}\label{e2equiv}
|j\omega-(1/2\pi{\rm i})\ln a-m|> cj^{-\nu}\qquad \mbox{ for all }j\in{\mathbb N}, m\in\mathbb Z,
\end{equation}
where $c$ and $\nu$ are some positive constants, then a formal power series solution of this linear $q$-difference equation converges \cite{b1}.
For the {\it nonlinear} equation \eqref{e1} an analogous sufficient condition \eqref{e2} of the convergence of its formal power series solution $\theta=\sum_{j=0}^{\infty}c_jz^j$ was obtained by Di Vizio \cite{dv}, where by $Q(\xi)$ one means a non-zero polynomial (depending on $\theta$) of degree
$\leqslant n$ whose coefficient $A_k$ is determined by the equality
$$
F'_{y_k}(z,\theta,\sigma\theta,\ldots,\sigma^n\theta)=A_kz^d+B_kz^{d+1}+\ldots, \qquad k=0,1,\ldots,n,
$$
which generalizes \eqref{lin}.

Now we pass to {\it generalized} formal power series solutions \eqref{eq1} of the equation \eqref{e1} whose convergence we continue to study in this paper.
In the situations where the small divisors phenomenon does not arise the convergence of such series has already been studied enough: this is again the case of $|q|\neq 1$ for solutions \eqref{eq1} with {\it real} power exponents $\lambda_j$'s (see a recent work by Barbe, Cano, Fortuny Ayuso, McCormick \cite{bcfm}) and a more general case of solutions \eqref{eq1} with complex power exponents $\lambda_j$'s where all the $q^{\lambda_j}$'s lie {\it strictly inside} or {\it strictly outside} the unit circle (or equivalently, where all the $\lambda_j$'s lie {\it strictly above} or {\it strictly under} the line passing through $0\in\mathbb C$ and having the slope $\ln |q|/\arg q$), which was considered in our previous paper \cite{GGL}. Note that the last condition of the placement of all the $q^{\lambda_j}$'s with respect to the unit circle (or of all the $\lambda_j$'s with respect to the line having the slope $\ln |q|/\arg q$) is invariant under the choice of the value of $\ln q$, since fixing another value of $\arg q$ one naturally comes to a different sequence of power exponents $\lambda_j$ according to a general approach to the construction of generalized power series solutions of \eqref{e1}.

Here we study the convergence of the formal solution \eqref{eq1} of the equation \eqref{e1} in the most general situation, that is, with the $\lambda_j$'s arbitrary placed. This is the situation where the small divisor phenomenon may arise. Assume that each $F'_{y_k}(z,\varphi,\sigma\varphi,\ldots,\sigma^n\varphi)$ is of the form
$$
\frac{\partial F}{\partial y_k}(z,\varphi,\sigma\varphi,\ldots,\sigma^n\varphi)=A_kz^{\alpha}+B_kz^{\alpha_k}+\ldots, \qquad {\rm Re}\,\alpha_k>
{\rm Re}\,\alpha\geqslant0,
$$
$\alpha\in\mathbb C$ being the same for all $k=0,1,\ldots,n$, and at least one of the $A_k$'s being non-zero. Then under a generic assumption on the power exponents $\lambda_j$ of \eqref{eq1} that, starting with some $j_0\in{\mathbb Z}_+$, the $q^{\lambda_j}$'s are not the roots of a non-zero polynomial
$$
L(\xi)=A_n\xi^n+\ldots+A_1\xi+A_0
$$
(depending on $\varphi$) of degree $\leqslant n$, one can assert that all $\lambda_j-\lambda_{j_0}$, $j>j_0$, belong to a finitely generated additive semi-group $\Gamma\subset\mathbb C$ whose generators $\alpha_1,\ldots,\alpha_s$, say, all have a positive real part (see Lemmas 1, 2 in \cite{GGL}). Thus we may initially consider the formal solution \eqref{eq1} in the form
\begin{equation}\label{gensol}
\varphi=\sum_{j=0}^{\infty}c_jz^{\lambda_j}=\sum_{j=0}^{j_0}c_jz^{\lambda_j}+\sum_{(m_1,\ldots,m_s)\in{\mathbb Z}_+^s\setminus\{0\}} c_{m_1,\ldots,m_s}z^{\lambda_{j_0}+m_1\alpha_1+\ldots+m_s\alpha_s}.
\end{equation}
Then we have the following theorem on the convergence of $\varphi$, which is the main result of the present paper.
\begin{theorem}\label{th1}
Let the generalized formal power series \eqref{gensol} satisfy \eqref{e1}. If $\deg{L}=n$, $L(0)\neq0$, and for each root $\xi=a$ of the polynomial $(\xi-q^{\lambda_{j_0}})L(\xi)$ the following diophantine condition is fulfilled:
\begin{equation}\label{smalldiv+}
|(\lambda_{j_0}+m_1\alpha_1+\ldots+m_s\alpha_s)\ln q-\ln a-2\pi m{\rm i}|>c\,|m_1+\ldots+m_s|^{-\nu} \quad \mbox{\rm for all}\; m_i\in{\mathbb Z}_+,\; 
m\in{\mathbb Z}
\end{equation}
$($with the exception of $m_1=\ldots=m_s=0)$, where $c$ and $\nu$ are some positive constants, then \eqref{gensol} has a non-zero radius of convergence $($that is, it converges uniformly in any sector $S\subset\mathbb C$ of sufficiently small radius with the vertex at the origin and of the opening less than $2\pi$ defining there a germ of a holomorphic function$)$.
\end{theorem}

The diophantine condition of Theorem \ref{th1} is generically fulfilled. As for concrete examples, one can apply in particular Schmidt's result \cite{Schmidt2} from which it follows that (\ref{smalldiv+}) holds for $a=q^{\lambda_{j_0}}$, if
\begin{enumerate}
\item the real parts of all $\frac1{2\pi\rm i}\alpha_1\ln q,\ldots,\frac1{2\pi\rm i}\alpha_s\ln q$ are algebraic and together with $1$ linearly independent over $\mathbb Z$ or
\item the imaginary parts of all $\frac1{2\pi\rm i}\alpha_1\ln q,\ldots,\frac1{2\pi\rm i}\alpha_s\ln q$ are algebraic and linearly independent over $\mathbb Z$.
\end{enumerate}
(If $L$ has roots $\xi=a$ other than $q^{\lambda_{j_0}}$ then the number $\frac1{2\pi\rm i}\ln(q^{\lambda_{j_0}}/a)$ should be added to the set of numbers in the above conditions 1, 2 for sufficiency of (\ref{smalldiv+}) for each such $a\ne q^{\lambda_{j_0}}$.)

\section{Proof of Theorem 1}

According to \eqref{gensol} the generalized formal power series $\varphi$ is presented in the form
$$
\varphi=\varphi_0+z^{\lambda}\psi, \qquad \lambda:=\lambda_{j_0},
$$
where
$$
\varphi_0=\sum_{j=0}^{j_0}c_jz^{\lambda_j}\quad \mbox{and}\quad
\psi=\sum_{(m_1,\ldots,m_s)\in{\mathbb Z}_+^s\setminus\{0\}}c_{m_1,\ldots,m_s}z^{m_1\alpha_1+\ldots+m_s\alpha_s},
$$
furthermore by Lemma 1 of \cite{GGL} the generalized formal power series $\psi$ satisfies an equality
\begin{equation}\label{e7}
L(q^{\lambda}\sigma)\psi=M(z,\psi,\sigma\psi,\ldots,\sigma^n\psi),
\end{equation}
where $M(z,u_0,u_1,\ldots,u_n)$ is a finite linear combination of monomials of the form
$$
z^{\alpha}u_0^{p_0}u_1^{p_1}\ldots u_n^{p_n}, \qquad \alpha\in\Gamma, \quad p_i\in\Z_+.
$$
Thus to prove Theorem \ref{th1}, it is sufficient to prove the convergence of $\psi$. 

Denoting $q_i=q^{\alpha_i}$ we rewrite the assumption \eqref{smalldiv+} of Theorem \ref{th1} in a multiplicative form as follows.

\begin{lemma}\label{l0}
The assumption \eqref{smalldiv+} implies the following two conditions for each root $\xi=a$ of the polynomial $L(\xi):$
\begin{eqnarray}
|q^{\lambda}q_1^{m_1}\ldots q_s^{m_s}-a| & > & 2^{-\nu}\,|m_1+\ldots+m_s|^{-\nu} \quad \mbox{\rm for all}\; m_i\in{\mathbb Z}_+, \label{smalldiv1} \\
|a\,q_1^{m_1}\ldots q_s^{m_s}-q^{\lambda}| & > & 2^{-\nu}\,|m_1+\ldots+m_s|^{-\nu} \quad \mbox{\rm for all}\; m_i\in{\mathbb Z}_-, \label{smalldiv2}
\end{eqnarray}
as well as the condition
\begin{eqnarray}
|q_1^{m_1}\ldots q_s^{m_s}-1|\;>\;2^{-\nu}\,|m_1+\ldots+m_s|^{-\nu} & & \mbox{\rm for all}\; m_i\in{\mathbb Z}_+ \; \mbox{\rm or} \label{smalldiv3} \\
                                                                    & & \mbox{\rm for all}\; m_i\in{\mathbb Z}_- \nonumber
\end{eqnarray}
$($with the exception of $m_1=\ldots=m_s=0)$, where $\nu$ is some positive constant $($maybe bigger than that in \eqref{smalldiv+}$)$.
\end{lemma}

\noindent{\bf Proof.} We will use the inequality $|e^z-1|>C|z|$ for $|z|<\varepsilon$ small enough.

1) Let $m_1,\ldots,m_s\in{\mathbb Z}_+$. If they are such that
$$
\Bigl|\sum_{i=1}^s m_i\alpha_i\ln q+\lambda\ln q-\ln a-2\pi m{\rm i}\Bigr|<\varepsilon
$$
for some $m\in\mathbb Z$, then we obtain
$$
|q^{\lambda}q_1^{m_1}\ldots q_s^{m_s}-a|=|a|\cdot|q^{\lambda}q_1^{m_1}\ldots q_s^{m_s}a^{-1}-1|=|a|\cdot|e^{\sum_{i=1}^s m_i\alpha_i\ln q+\lambda\ln q-\ln a-2\pi m{\rm i}}-1|>
$$
$$
>|a|\cdot C\Bigl|\sum_{i=1}^s m_i\alpha_i\ln q+\lambda\ln q-\ln a-2\pi m{\rm i}\Bigr|>c_1\,|m_1+\ldots+m_s|^{-\nu},
$$
whence \eqref{smalldiv1} follows. In the case where $m_1,\ldots,m_s$ are such that
$$
\Bigl|\sum_{i=1}^s m_i\alpha_i\ln q+\lambda\ln q-\ln a-2\pi m{\rm i}\Bigr|\geqslant\varepsilon
$$
for every $m\in\mathbb Z$, we obtain
$$
|q^{\lambda}q_1^{m_1}\ldots q_s^{m_s}-a|=|a|\cdot|e^{\sum_{i=1}^s m_i\alpha_i\ln q+\lambda\ln q-\ln a}-1|>c_2>0,
$$
whence \eqref{smalldiv1} follows again.

2) Let $m_1,\ldots,m_s\in{\mathbb Z}_-$. If they are such that
$$
\Bigl|\sum_{i=1}^s m_i\alpha_i\ln q-\lambda\ln q+\ln a-2\pi m{\rm i}\Bigr|<\varepsilon
$$
for some $m\in\mathbb Z$, then we obtain
$$
|a\,q_1^{m_1}\ldots q_s^{m_s}-q^{\lambda}|=|q^{\lambda}|\cdot|a\,q_1^{m_1}\ldots q_s^{m_s}q^{-\lambda}-1|=|q^{\lambda}|\cdot |e^{\sum_{i=1}^s m_i\alpha_i\ln q-\lambda\ln q+\ln a-2\pi m{\rm i}}-1|>
$$
$$
>|q^{\lambda}|\cdot C\Bigl|\sum_{i=1}^s m_i\alpha_i\ln q-\lambda\ln q+\ln a-2\pi m{\rm i}\Bigr|=|q^{\lambda}|\cdot C\Bigl|\sum_{i=1}^s-m_i\alpha_i\ln q+\lambda\ln q -\ln a+2\pi m{\rm i}\Bigr|>
$$
$$
>\tilde c_1\,|-m_1-\ldots-m_s|^{-\nu},
$$
whence \eqref{smalldiv2} follows. In the case where $m_1,\ldots,m_s$ are such that
$$
\Bigl|\sum_{i=1}^s m_i\alpha_i\ln q-\lambda\ln q+\ln a-2\pi m{\rm i}\Bigr|\geqslant\varepsilon
$$
for every $m\in\mathbb Z$, we obtain
$$
|a\,q_1^{m_1}\ldots q_s^{m_s}-q^{\lambda}|=|q^{\lambda}|\cdot|e^{\sum_{i=1}^s m_i\alpha_i\ln q-\lambda\ln q+\ln a}-1|>\tilde c_2>0,
$$
whence \eqref{smalldiv2} follows again.

3) The estimate \eqref{smalldiv3} is obtained in an absolutely similar manner if we use the condition
$$
|(m_1\alpha_1+\ldots+m_s\alpha_s)\ln q-2\pi m{\rm i}|>c\,|m_1+\ldots+m_s|^{-\nu} \quad \mbox{\rm for all}\; m_i\in{\mathbb Z}_+,\; m\in{\mathbb Z},
$$
coming from (\ref{smalldiv+}) with $a=q^{\lambda_{j_0}}$. The lemma is proved.

\subsection{A reduced equation and its majorant one}

We will always use boldface letters for multi-indices of length $s$, such as ${\bf k}=(k_1,\ldots,k_s)$,
${\bf m}=(m_1,\ldots,m_s)$, etc. As usually, $|{\bf k}|, |{\bf m}|$, etc., will denote the sum of the coordinates of multi-indices.

As we know from \cite{GGL}, the representation
\begin{equation}\label{formsol}
\tilde\psi=\sum_{|{\bf m}|>0}c_{\bf m}\,z_1^{m_1}\ldots z_s^{m_s}
\end{equation}
of the formal series $\psi$ by an $s$-variate formal Taylor series satisfies the equality (the {\it reduced} equation)
\begin{equation}\label{a}
L(q^{\lambda}\tilde\sigma)\tilde\psi=\widetilde M(z_1,\ldots,z_s,\tilde\psi,\tilde\sigma\tilde\psi,\ldots,
\tilde\sigma^n\tilde\psi),
\end{equation}
where
\begin{eqnarray}
\tilde\sigma^j\tilde\psi&=&\sum_{|{\bf m}|>0}q_1^{jm_1}\ldots q_s^{jm_s}c_{\bf m}\,z_1^{m_1}\ldots z_s^{m_s}, \quad j=0,1,\ldots,n, \label{b} \\
L(q^{\lambda}\tilde\sigma)\tilde\psi&=&\sum_{|{\bf m}|>0}L(q^{\lambda}q_1^{m_1}\ldots q_s^{m_s})c_{\bf m}\, z_1^{m_1}\ldots z_s^{m_s}, \label{c}
\end{eqnarray}
and $\widetilde M(z_1,\ldots,z_s,u_0,\ldots,u_n)$ is a polynomial of the form
\begin{equation}\label{d}
\widetilde M(z_1,\ldots,z_s,u_0,\ldots,u_n)=\sum_{|{\bf k}|>0,\,p_i\geqslant0}A_{{\bf k},\,p_0,\ldots,p_n}\,z_1^{k_1}
\ldots z_s^{k_s}u_0^{p_0}\ldots u_n^{p_n}.
\end{equation}
Note that $L(q^{\lambda}q_1^{m_1}\ldots q_s^{m_s})\ne0$ for every $\bf m$ by \eqref{smalldiv1}.

To prove the convergence of $\tilde\psi$ in some neighbourhood of the origin, we construct an equation
\begin{eqnarray}\label{majorant}
\tilde\nu\,W=\sum_{|{\bf k}|>0,\,p_i\geqslant0}|A_{{\bf k},\,p_0,\ldots,p_n}|\,z_1^{k_1}\ldots z_s^{k_s}
W^{p_0}\ldots W^{p_n},
\end{eqnarray}
whose right-hand side is obtained from the polynomial $\widetilde M$ by the change of its coefficients $A_{{\bf k},\,p_0,\ldots,p_n}$ to their absolute values and of all the $u_j$'s to the one variable $W$. The number $0<\tilde\nu\leqslant1$ will be defined further.

\subsection{The majorant solution}

Due to the implicit function theorem equation (\ref{majorant}) possesses a unique solution $W=\Psi(z_1,\ldots,z_s)$ holomorphic near the origin,
\begin{equation}\label{solmaj}
\Psi=\sum\limits_{|{\bf m}|>0}C_{\bf m}\,z_1^{m_1}\ldots z_s^{m_s}, \qquad \Psi(0,\ldots,0)=0,
\end{equation}
which, as we will prove further in this and the next sections, is majorant for $\tilde\psi$. By this we mean that the ratio $|c_{\bf m}|/C_{\bf m}$ has an exponential growth at most with respect to $|{\bf m}|$ (the coefficients $C_{\bf m}$ are real non-negative numbers, by construction). 

By notation ${\bf m}>{\bf p}$ we will mean that $m_i\geqslant p_i$ for all $i$, and $|{\bf m}|>|{\bf p}|$. Introducing
$$
\epsilon_{\bf m}=|L(q^{\lambda}q_1^{m_1}\ldots q_s^{m_s})|^{-1}, \qquad s_{\bf m}=\max(1, |q_1|^{m_1}\ldots|q_s|^{m_s}),
$$
consider also the sequence of positive numbers $\delta_{\bf m}$ recurrently defined in the following way. For $|{\bf m}|=1$ put $\delta_{\bf m}=1$ and for every
$|{\bf m}|>1$ let $\mu_{\bf m}$ denote the maximum of all products
\begin{eqnarray*}
s^n_{{\bf m}^{(1)}}\delta_{{\bf m}^{(1)}}\ldots s^n_{{\bf m}^{(k)}}\delta_{{\bf m}^{(k)}}, & \mbox{with} & {{\bf m}^{(1)}}+\ldots+{{\bf m}^{(k)}}={\bf m},
\;k\geqslant2,\\
                                                                                           &           & {{\bf m}^{(1)}}<{\bf m},\ldots,{{\bf m}^{(k)}}<{\bf m}.
\end{eqnarray*}
Then define
$$
\delta_{\bf m}=\epsilon_{\bf m}\mu_{\bf m}.
$$

Now let us choose a number $\tilde\nu\leqslant1$ satisfying the condition
$$
\tilde\nu\,\max_{|{\bf m}|=1}\epsilon_{\bf m}\leqslant1
$$
and after such choice prove the following lemma.

\begin{lemma}\label{lemmamaj}
The estimate
$$
|c_{\bf m}|\leqslant\delta_{\bf m}C_{\bf m}
$$
holds for all the coefficients of the series \eqref{formsol}, \eqref{solmaj}.
\end{lemma}

\noindent{\bf Proof.} The estimate is true in the case $|{\bf m}|=1$:
$$
\epsilon_{\bf m}^{-1}|c_{\bf m}|=|A_{{\bf m},0,\ldots,0}|=\tilde\nu\,C_{\bf m}\quad\Longrightarrow\quad
|c_{\bf m}|=\epsilon_{\bf m}\tilde\nu\,C_{\bf m}\leqslant\delta_{\bf m}C_{\bf m}.
$$
Further we assume $|{\bf m}|>1$ and apply induction.

Using the notations
$$
\frac{\partial^{\bf m}}{{\bf m}!}=\frac1{m_1!\ldots m_s!}\frac{\partial^{m_1}}{\partial z_1^{m_1}}\ldots
\frac{\partial^{m_s}}{\partial z_s^{m_s}}, \qquad {\bf q}^{\bf m}=q_1^{m_1}\ldots q_s^{m_s},
$$
and defining
$$
\phi=\widetilde M(z_1,\ldots,z_s,\tilde\psi,\tilde\sigma\tilde\psi,\ldots,\tilde\sigma^n\tilde\psi),
$$
one has in view of \eqref{a}, \eqref{c}, \eqref{d}:
$$
L(q^{\lambda}\,{\bf q}^{\bf m})c_{\bf m}=\Bigl.\frac{\partial^{\bf m}\phi}{{\bf m}!}\Bigr|_{z_1=\ldots= z_s=0}=
$$
\begin{equation}\label{relation1}
=\sum_{0<{\bf k}\leqslant{\bf m},\,p_i\geqslant0}A_{{\bf k},\,p_0,\ldots,p_n}\sum_{{\bf l}^{(0)}+\ldots+{\bf l}^{(n)}=
{\bf m}-{\bf k}}\frac{\partial^{{\bf l}^{(0)}}\tilde\psi^{p_0}}{{\bf l}^{(0)}!}\ldots
\Bigl.\frac{\partial^{{\bf l}^{(n)}}(\tilde\sigma^n\tilde\psi)^{p_n}}{{\bf l}^{(n)}!}\Bigr|_{z_1=\ldots= z_s=0}.
\end{equation}
On the other hand, defining
$$
\Phi=\sum_{|{\bf k}|>0,\,p_i\geqslant0}|A_{{\bf k},\,p_0,\ldots,p_n}|\,z_1^{k_1}\ldots z_s^{k_s}\Psi^{p_0}\Psi^{p_1}\ldots\Psi^{p_n},
$$
in view of \eqref{majorant} one has for the coefficients of the series \eqref{solmaj} the equality
$$
\tilde\nu\,C_{\bf m}=\Bigl.\frac{\partial^{\bf m}\Phi}{{\bf m}!}\Bigr|_{z_1=\ldots= z_s=0}=
$$
\begin{equation}\label{relation2}
=\sum_{0<{\bf k}\leqslant{\bf m},\,p_i\geqslant0}|A_{{\bf k},\,p_0,\ldots,p_n}|\sum_{{\bf l}^{(0)}+\ldots+{\bf l}^{(n)}=
{\bf m}-{\bf k}}\frac{\partial^{{\bf l}^{(0)}}\Psi^{p_0}}{{\bf l}^{(0)}!}\ldots
\Bigl.\frac{\partial^{{\bf l}^{(n)}}\Psi^{p_n}}{{\bf l}^{(n)}!}\Bigr|_{z_1=\ldots= z_s=0}.
\end{equation}

Noting that each
$$
\Bigl.\frac{\partial^{{\bf l}^{(j)}}(\tilde\sigma^j\tilde\psi)^{p_j}}{{\bf l}^{(j)}!}\Bigr|_{z_1=\ldots= z_s=0}=
\sum_{\boldsymbol{\lambda}^{(1)}+\ldots+\boldsymbol{\lambda}^{(p_j)}={\bf l}^{(j)}}
{\bf q}^{j\boldsymbol{\lambda}^{(1)}}c_{\boldsymbol{\lambda}^{(1)}}\ldots {\bf q}^{j\boldsymbol{\lambda}^{(p_j)}}c_{\boldsymbol{\lambda}^{(p_j)}}, \qquad j=0,1,\ldots,n,
$$
by its absolute value does not exceed $\sum\limits_{\boldsymbol{\lambda}^{(1)}+\ldots+{\boldsymbol{\lambda}}^{(p_j)}={\bf l}^{(j)}}
s_{\boldsymbol{\lambda}^{(1)}}^n|c_{{\boldsymbol{\lambda}}^{(1)}}|\ldots s_{\boldsymbol{\lambda}^{(p_j)}}^n|c_{{\boldsymbol{\lambda}}^{(p_j)}}|$, one applies induction and obtains
$$
\Bigl|\frac{\partial^{{\bf l}^{(j)}}(\tilde\sigma^j\tilde\psi)^{p_j}}{{\bf l}^{(j)}!}\Bigr|_{z_1=\ldots= z_s=0}\leqslant
\sum\limits_{\boldsymbol{\lambda}^{(1)}+\ldots+{\boldsymbol{\lambda}}^{(p_j)}={\bf l}^{(j)}}
s_{\boldsymbol{\lambda}^{(1)}}^n\delta_{\boldsymbol{\lambda}^{(1)}}\ldots s_{\boldsymbol{\lambda}^{(p_j)}}^n\delta_{\boldsymbol{\lambda}^{(p_j)}} C_{{\boldsymbol{\lambda}}^{(1)}}\ldots C_{{\boldsymbol{\lambda}}^{(p_j)}}.
$$
We have the following dichotomy: if $p_j=1$ then the last sum is equal to
$$
s_{{\bf l}^{(j)}}^n\,\delta_{{\bf l}^{(j)}}C_{{\bf l}^{(j)}}=s_{{\bf l}^{(j)}}^n\,\delta_{{\bf l}^{(j)}}
\Bigl.\frac{\partial^{{\bf l}^{(j)}}\Psi^{p_j}}{{\bf l}^{(j)}!}\Bigr|_{z_1=\ldots=z_s=0},
$$
whereas in the case $p_j\geqslant2$, by the definition of the sequence $\mu_{\bf m}$, that sum does not exceed
$$
\mu_{{\bf l}^{(j)}}\sum\limits_{\boldsymbol{\lambda}^{(1)}+\ldots+\boldsymbol{\lambda}^{(p_j)}={\bf l}^{(j)}}C_{\boldsymbol{\lambda}^{(1)}}\ldots
C_{\boldsymbol{\lambda}^{(p_j)}}=\mu_{{\bf l}^{(j)}}\Bigl.\frac{\partial^{{\bf l}^{(j)}}\Psi^{p_j}}{{\bf l}^{(j)}!}\Bigr|_{z_1=\ldots=z_s=0}.
$$
Thus we come to the estimate
$$
\Bigl|\frac{\partial^{{\bf l}^{(0)}}\tilde\psi^{p_0}}{{\bf l}^{(0)}!}\ldots
\frac{\partial^{{\bf l}^{(n)}}(\tilde\sigma^n\tilde\psi)^{p_n}}{{\bf l}^{(n)}!}\Bigr|_{z_1=\ldots= z_s=0}\leqslant\mu_{\bf m-k}\,
\frac{\partial^{{\bf l}^{(0)}}\Psi^{p_0}}{{\bf l}^{(0)}!}\ldots\Bigl.\frac{\partial^{{\bf l}^{(n)}}\Psi^{p_n}}{{\bf l}^{(n)}!}\Bigr|_{z_1=\ldots= z_s=0}
$$
for each set of multi-indices ${\bf l}^{(0)},{\bf l}^{(1)},\ldots,{\bf l}^{(n)}$ with ${\bf l}^{(0)}+\ldots+{\bf l}^{(n)}={\bf m-k}$. From this estimate, in view of
\eqref{relation1}, it follows that
$$
\epsilon_{\bf m}^{-1}|c_{\bf m}|\leqslant\sum_{0<{\bf k}\leqslant{\bf m},\,p_i\geqslant0}|A_{{\bf k},\,p_0,\ldots,p_n}|
\sum_{{\bf l}^{(0)}+\ldots+{\bf l}^{(n)}={\bf m}-{\bf k}}\Bigl|\frac{\partial^{{\bf l}^{(0)}}\tilde\psi^{p_0}}{{\bf l}^{(0)}!}\ldots
\frac{\partial^{{\bf l}^{(n)}}(\tilde\sigma^n\tilde\psi)^{p_n}}{{\bf l}^{(n)}!}\Bigr|_{z_1=\ldots= z_s=0}\leqslant
$$
\begin{equation}\label{est3}
\leqslant\sum_{0<{\bf k}\leqslant{\bf m},\,p_i\geqslant0}|A_{{\bf k},\,p_0,\ldots,p_n}|\,\mu_{\bf m-k}\sum_{{\bf l}^{(0)}+\ldots+{\bf l}^{(n)}={\bf m}-{\bf k}}
\frac{\partial^{{\bf l}^{(0)}}\Psi^{p_0}}{{\bf l}^{(0)}!}\ldots\Bigl.\frac{\partial^{{\bf l}^{(n)}}\Psi^{p_n}}{{\bf l}^{(n)}!}\Bigr|_{z_1=\ldots= z_s=0}.
\end{equation}
Note that $\mu_{{\bf m}-{\bf k}}\leqslant \mu_{\bf m}$: for $|{\bf k}|=1$ this follows from the inequalities
$$
\mu_{\bf m-k}\leqslant s_{\bf k}^n\,\mu_{\bf m-k}=s_{\bf k}^n\,\delta_{\bf k}\,\mu_{\bf m-k}\leqslant\mu_{\bf m}
$$
($\delta_{\bf k}=1$ in this case, $s_{\bf k}\geqslant 1$), for any other $\bf k$ one may proceed recursively. Therefore, by \eqref{est3} and \eqref{relation2} we finally obtain
\begin{eqnarray*}
\epsilon_{\bf m}^{-1}|c_{\bf m}|&\leqslant&\mu_{\bf m}\sum_{0<{\bf k}\leqslant{\bf m},\,p_i\geqslant0}|A_{{\bf k},\,p_0,\ldots,p_n}|
\sum_{{\bf l}^{(0)}+\ldots+{\bf l}^{(n)}={\bf m}-{\bf k}}
\frac{\partial^{{\bf l}^{(0)}}\Psi^{p_0}}{{\bf l}^{(0)}!}\ldots\Bigl.\frac{\partial^{{\bf l}^{(n)}}\Psi^{p_n}}{{\bf l}^{(n)}!}\Bigr|_{z_1=\ldots= z_s=0}
=\\ &=&\mu_{\bf m}\,\tilde\nu\,C_{\bf m}\leqslant\mu_{\bf m}C_{\bf m},
\end{eqnarray*}
whence the required estimate of the lemma follows.
\medskip

Now in view of Lemma \ref{lemmamaj}, to prove the convergence of the $s$-variate power series (\ref{formsol}) representing the generalized power series $\psi$ which satisfies (\ref{e7}), it is sufficient to establish that the sequence $\delta_{\bf m}$ has the exponential growth at most. This will be precisely established in Lemma \ref{l3} preceded by several auxiliary estimates in spirit of Siegel.

\subsection{Preliminary estimates}

Assuming the polynomial $L$ to be monic, one decomposes $L(q^{\lambda}x)$ in the product of factors $q^{\lambda}x-a$, the $a$'s being the roots of $L$. Define
$$
\varepsilon_{\bf m}(a)=\bigl|q^{\lambda}q_1^{m_1}\ldots q_s^{m_s}-a\bigr|^{-1}.
$$

\begin{lemma}\label{l1}
The following estimates hold for each non-zero root $a$ of the polynomial $L$ and any ${\bf m}>{\bf p}:$
$$
s_{\bf m}\,\min(\varepsilon_{\bf m}(a),\varepsilon_{\bf p}(a))<2^{\nu+1}(|{\bf m}|-|{\bf p}|)^{\nu},
$$
$$
s_{\bf p}\,\min(\varepsilon_{\bf m}(a),\varepsilon_{\bf p}(a))<2^{\nu+1}(|{\bf m}|-|{\bf p}|)^{\nu},
$$
where $\nu$ is some positive constant $($maybe bigger than that in Lemma \ref{l0}$)$.
\end{lemma}

\noindent{\bf Proof.} For any ${\bf m}>{\bf p}$ we have the identities
\begin{equation}\label{id1}
a\,\bigl(q_1^{m_1-p_1}\ldots q_s^{m_s-p_s}-1\bigr)=\bigl(q^{\lambda}q_1^{m_1}\ldots q_s^{m_s}-a\bigr)-
\bigl(q^{\lambda}q_1^{p_1}\ldots q_s^{p_s}-a\bigr)q_1^{m_1-p_1}\ldots q_s^{m_s-p_s},
\end{equation}
\begin{equation}\label{id2}
a\,\bigl(q_1^{p_1-m_1}\ldots q_s^{p_s-m_s}-1\bigr)=\bigl(q^{\lambda}q_1^{p_1}\ldots q_s^{p_s}-a\bigr)-
\bigl(q^{\lambda}q_1^{m_1}\ldots q_s^{m_s}-a\bigr)q_1^{p_1-m_1}\ldots q_s^{p_s-m_s}.
\end{equation}
Therefore, if $|q_1^{m_1-p_1}\ldots q_s^{m_s-p_s}|\leqslant1$, then \eqref{id1} implies
$$
|a|\,|q_1^{m_1-p_1}\ldots q_s^{m_s-p_s}-1|\leqslant\varepsilon_{\bf m}^{-1}(a)+\varepsilon_{\bf p}^{-1}(a)\leqslant
\frac2{\min(\varepsilon_{\bf m}(a),\varepsilon_{\bf p}(a))}.
$$
On the other hand, if $|q_1^{m_1-p_1}\ldots q_s^{m_s-p_s}|\geqslant1$, then \eqref{id2} implies
$$
|a|\,|q_1^{p_1-m_1}\ldots q_s^{p_s-m_s}-1|\leqslant\varepsilon_{\bf p}^{-1}(a)+\varepsilon_{\bf m}^{-1}(a)\leqslant
\frac2{\min(\varepsilon_{\bf m}(a),\varepsilon_{\bf p}(a))}.
$$
In any case applying \eqref{smalldiv3}, for a non-zero root $a$ of $L$ one has
\begin{equation}\label{est1}
\min(\varepsilon_{\bf m}(a),\varepsilon_{\bf p}(a))<\frac{2^{\nu+1}}{|a|}\,(|{\bf m}|-|{\bf p}|)^{\nu}.
\end{equation}

For any ${\bf m}>{\bf p}$ we also have the identities
$$
\bigl(q^{\lambda}q_1^{m_1}\ldots q_s^{m_s}-a\bigr)-\bigl(q^{\lambda}q_1^{p_1}\ldots q_s^{p_s}-a\bigr)=
q^{\lambda}q_1^{p_1}\ldots q_s^{p_s}(q_1^{m_1-p_1}\ldots q_s^{m_s-p_s}-1),
$$
$$\bigl(q^{\lambda}q_1^{m_1}\ldots q_s^{m_s}-a\bigr)-\bigl(q^{\lambda}q_1^{p_1}\ldots q_s^{p_s}-a\bigr)=
q^{\lambda}q_1^{m_1}\ldots q_s^{m_s}(1-q_1^{p_1-m_1}\ldots q_s^{p_s-m_s}),
$$
which imply
$$
|q^\lambda|\,|q_1|^{p_1}\ldots|q_s|^{p_s}\,|q_1^{m_1-p_1}\ldots q_s^{m_s-p_s}-1|\leqslant\varepsilon_{\bf m}^{-1}(a)+\varepsilon_{\bf p}^{-1}(a)\leqslant
\frac2{\min(\varepsilon_{\bf m}(a),\varepsilon_{\bf p}(a))},
$$
$$
|q^{\lambda}|\,|q_1|^{m_1}\ldots|q_s|^{m_s}\,|q_1^{p_1-m_1}\ldots q_s^{p_s-m_s}-1|\leqslant\varepsilon_{\bf m}^{-1}(a)+\varepsilon_{\bf p}^{-1}(a)\leqslant
\frac2{\min(\varepsilon_{\bf m}(a),\varepsilon_{\bf p}(a))}.
$$
Therefore by \eqref{smalldiv3}, for any non-zero root $a$ of $L$ one has
\begin{equation}\label{est2}
|q_1|^{m_1}\ldots|q_s|^{m_s}\,\min(\varepsilon_{\bf m}(a),\varepsilon_{\bf p}(a))<\frac{2^{\nu+1}}{|q^{\lambda}|}\,(|{\bf m}|-|{\bf p}|)^{\nu},
\end{equation}
\begin{equation}\label{est2'}
|q_1|^{p_1}\ldots|q_s|^{p_s}\,\min(\varepsilon_{\bf m}(a),\varepsilon_{\bf p}(a))<\frac{2^{\nu+1}}{|q^{\lambda}|}\,(|{\bf m}|-|{\bf p}|)^{\nu}.
\end{equation}
Increasing if necessary the value of $\nu$, we come from (\ref{est1}), (\ref{est2}) and (\ref{est2'}) to the estimate of the lemma.
\medskip

\noindent{\bf Remark 1.} If all $|q_i|\geqslant1$, the assumption $a\ne0$ is unnecessary for the estimates of Lemmas \ref{l0} and \ref{l1} to be held.
\medskip

Further by $\nu$ we will always mean the constant of Lemma \ref{l1}. The latter implies the following analogue of Siegel's Lemma 2 from \cite{Siegel}.
\begin{lemma}\label{l2}
Let ${\bf m}^{(0)}>{\bf m}^{(1)}>\ldots>{\bf m}^{(r)}>0$, $r\geqslant0$. For any non-zero root $a$ of the polynomial $L$ there holds
$$
\prod\limits_{i=0}^rs_{{\bf m}^{(i)}}\varepsilon_{{\bf m}^{(i)}}(a)<N_1^{r+1}|{\bf m}^{(0)}|^{\nu}\prod\limits_{i=1}^r
(|{\bf m}^{(i-1)}|-|{\bf m}^{(i)}|)^{\nu}, \qquad N_1=2^{2\nu+1}.
$$
\end{lemma}

\noindent{\bf Proof.} The assertion is true in the case $r=0$ in view of (\ref{smalldiv1}), (\ref{smalldiv2}); assume $r>0$ and apply induction. Let
$\varepsilon_{{\bf m}^{(i)}}(a)$, $i=0,1,\ldots,r$, have its minimum value for $i=h$. Then Lemma \ref{l1} yields
\begin{equation}\label{min}
s_{{\bf m}^{(h)}}\varepsilon_{{\bf m}^{(h)}}(a)<2^{\nu+1}\min\bigl\{(|{\bf m}^{(h-1)}|-|{\bf m}^{(h)}|)^{\nu},
(|{\bf m}^{(h)}|-|{\bf m}^{(h+1)}|)^{\nu}\bigr\},
\end{equation}
if we define moreover $|{\bf m}^{(-1)}|=\infty$ and $|{\bf m}^{(r+1)}|=-\infty$. On the other hand, the lemma being true for $r-1$ instead of $r$, we have
$$
s_{{\bf m}^{(h)}}^{-1}\varepsilon_{{\bf m}^{(h)}}^{-1}(a)\prod\limits_{i=0}^{r}s_{{\bf m}^{(i)}}\varepsilon_{{\bf m}^{(i)}}(a)<
N_1^r|{\bf m}^{(0)}|^{\nu}\frac{(|{\bf m}^{(h-1)}|-|{\bf m}^{(h+1)}|)^{\nu}\prod\limits_{i=1}^r(|{\bf m}^{(i-1)}|-|{\bf m}^{(i)}|)^{\nu}}
{(|{\bf m}^{(h-1)}|-|{\bf m}^{(h)}|)^{\nu}(|{\bf m}^{(h)}|-|{\bf m}^{(h+1)}|)^{\nu}}=
$$
$$
=N_1^r|{\bf m}^{(0)}|^{\nu}\Bigl(\frac{1}{|{\bf m}^{(h-1)}|-|{\bf m}^{(h)}|}+\frac{1}{|{\bf m}^{(h)}|-|{\bf m}^{(h+1)}|}\Bigr)^{\nu}
\prod\limits_{i=1}^r(|{\bf m}^{(i-1)}|-|{\bf m}^{(i)}|)^{\nu}\leqslant
$$
$$
\leqslant\frac{2^\nu N_1^r}{\min\bigl\{(|{\bf m}^{(h-1)}|-|{\bf m}^{(h)}|)^{\nu},(|{\bf m}^{(h)}|-|{\bf m}^{(h+1)}|)^{\nu}\bigr\}}\,|{\bf m}^{(0)}|^{\nu} \prod\limits_{i=1}^r(|{\bf m}^{(i-1)}|-|{\bf m}^{(i)}|)^{\nu},
$$
and the statement of the lemma follows from \eqref{min}.

\begin{corol}\label{cor1}
In the case where zero is not a root of the polynomial $L$ and $\deg L=n$, one has the following estimate:
$$
\prod\limits_{i=0}^rs^n_{{\bf m}^{(i)}}\epsilon_{{\bf m}^{(i)}}<
\Bigl(N_1^{r+1}|{\bf m}^{(0)}|^{\nu}\prod\limits_{i=1}^r(|{\bf m}^{(i-1)}|-|{\bf m}^{(i)}|)^{\nu}\Bigr)^n.
$$
\end{corol}

\noindent{\bf Proof.} Follows from the equality
$$
s^n_{{\bf m}^{(i)}}\epsilon_{{\bf m}^{(i)}}=s_{{\bf m}^{(i)}}\varepsilon_{{\bf m}^{(i)}}(a_1)\ldots s_{{\bf m}^{(i)}}\varepsilon_{{\bf m}^{(i)}}(a_n),
$$
where $a_1,\ldots,a_n\ne0$ are the roots of $L$.
\medskip

\noindent{\bf Remark 2.} If all $|q_i|\geqslant1$, like in Remark 1, the assumption $a\ne0$ is unnecessary for the estimate of Lemma \ref{l2} to be held. Hence Corollary \ref{cor1} is true in this case even if $L(0)=0$.

On the other hand, if all $|q_i|\leqslant1$ then Corollary \ref{cor1} is true for $d=\deg L<n$ as well,
since $s^n_{\bf m}=s^d_{\bf m}=1$ for any $\bf m$ in this case.

\subsection{The main estimate}

Now we finally prove that the sequence $\delta_{\bf m}$ has the exponential growth at most.
\begin{lemma}\label{l3}
\begin{equation}\label{lemma3}
s_{\bf m}^n\delta_{\bf m}\leqslant|{\bf m}|^{-2\nu n}N_2^{|{\bf m}|-1}Q^{|{\bf m}|},
\end{equation}
where $N_2=8^{\nu n}N_1^n$, $Q=\max(1,|q_1|^n,\ldots,|q_s|^n)$.
\end{lemma}

\noindent
{\bf Proof.} The assertion is true for $|{\bf m}|=1$, since $s_{\bf m}^n\leqslant Q$, $\delta_{\bf m}=1$ in this case. Assume $|{\bf m}|>1$ and apply induction.

The numbers $\alpha_{\bf m}=|{\bf m}|^{-2\nu n}N_2^{|{\bf m}|-1}Q^{|{\bf m}|}$ satisfy the inequalities
$$
\frac{\alpha_{\bf m}\alpha_{\bf l}}{\alpha_{{\bf m}+{\bf l}}}=(|{\bf m}|^{-1}+|{\bf l}|^{-1})^{2\nu n}N_2^{-1}\leqslant
2^{2\nu n}N_2^{-1}<1, \qquad |{\bf m}|\geqslant1, |{\bf l}|\geqslant1,
$$
and consequently
\begin{equation}\label{estimate}
s_{{\bf m}^{(1)}}^n\delta_{{\bf m}^{(1)}}\ldots s_{{\bf m}^{(f)}}^n\delta_{{\bf m}^{(f)}}\leqslant|{\bf j}|^{-2\nu n}N_2^{|{\bf j}|-1}Q^{|{\bf j}|} \qquad
({\bf m}^{(1)}+\ldots+{\bf m}^{(f)}={\bf j}, \; |{\bf j}|<|{\bf m}|).
\end{equation}
By definition of $\delta_{\bf m}$, there exists a decomposition
$$
s_{\bf m}^n\delta_{\bf m}=\bigl(s_{\bf m}^n\epsilon_{\bf m}\bigr)s_{{\bf g}^{(1)}}^n\delta_{{\bf g}^{(1)}}\ldots s_{{\bf g}^{(\alpha)}}^n\delta_{{\bf g}^{(\alpha)}}, \qquad \alpha\geqslant2,
$$
with ${\bf g}^{(1)}+\ldots+{\bf g}^{(\alpha)}={\bf m}$, and ${\bf g}^{(1)}<{\bf m},\ldots,{\bf g}^{(\alpha)}<{\bf m}$.
In the case $|{\bf g}^{(1)}|>|{\bf m}|/2$ (and hence $|{\bf g}^{(2)}|+\ldots+|{\bf g}^{(\alpha)}|<|{\bf m}|/2$) we use this formula with ${\bf g}^{(1)}$ instead of ${\bf m}$ and find a decomposition
$$
s_{{\bf g}^{(1)}}^n\delta_{{\bf g}^{(1)}}=\bigl(s_{{\bf g}^{(1)}}^n\epsilon_{{\bf g}^{(1)}}\bigr)s_{{\bf h}^{(1)}}^n\delta_{{\bf h}^{(1)}}\ldots
s_{{\bf h}^{(\beta)}}^n\delta_{{\bf h}^{(\beta)}}, \qquad \beta\geqslant 2,
$$
with ${\bf h}^{(1)}+\ldots+{\bf h}^{(\beta)}={\bf g}^{(1)}$, and ${\bf h}^{(1)}<{\bf g}^{(1)},\ldots,{\bf h}^{(\beta)}<{\bf g}^{(1)}$.
If also $|{\bf h}^{(1)}|>|{\bf m}|/2$ (and hence $|{\bf h}^{(2)}|+\ldots+|{\bf h}^{(\beta)}|<|{\bf m}|/2$) we decompose again
$$
s_{{\bf h}^{(1)}}^n\delta_{{\bf h}^{(1)}}=\bigl(s_{{\bf h}^{(1)}}^n\epsilon_{{\bf h}^{(1)}}\bigr)s_{{\bf i}^{(1)}}^n\delta_{{\bf i}^{(1)}}\ldots
s_{{\bf i}^{(\gamma)}}^n\delta_{{\bf i}^{(\gamma)}}, \qquad \gamma\geqslant 2,
$$
with ${\bf i}^{(1)}+\ldots+{\bf i}^{(\gamma)}={\bf h}^{(1)}$, and ${\bf i}^{(1)}<{\bf h}^{(1)},\ldots,{\bf i}^{(\gamma)}<{\bf h}^{(1)}$, and so on. Taking
$$
{\bf m}^{(0)}={\bf m},\; {\bf m}^{(1)}={\bf g}^{(1)},\; {\bf m}^{(2)}={\bf h}^{(1)},\;\ldots,\;{\bf m}^{(r)}=\dots,
$$
$$
{\bf m}^{(0)}>{\bf m}^{(1)}>{\bf m}^{(2)}>\ldots>{\bf m}^{(r)},
$$
${\bf m}^{(r)}$ being the last multi-index in this chain with $|{\bf m}^{(r)}|>|{\bf m}|/2$, we obtain in this manner the formula
\begin{equation}\label{decomp}
s_{\bf m}^n{\delta}_{\bf m}=\bigl(s_{{\bf m}^{(0)}}^n\epsilon_{{\bf m}^{(0)}}\Delta_0\bigr)\bigl(s_{{\bf m}^{(1)}}^n\epsilon_{{\bf m}^{(1)}}\Delta_1\bigr)\ldots
\bigl(s_{{\bf m}^{(r)}}^n\epsilon_{{\bf m}^{(r)}}\Delta_r\bigr),
\end{equation}
where
$$
\begin{array}{lllll}
\Delta_0&=&s_{{\bf g}^{(2)}}^n\delta_{{\bf g}^{(2)}}\ldots s_{{\bf g}^{(\alpha)}}^n\delta_{{\bf g}^{(\alpha)}}, & & {\bf g}^{(2)}+\ldots+{\bf g}^{(\alpha)}=
                                                                                                                    {\bf m}^{(0)}-{\bf m}^{(1)},\\
\Delta_1&=&s_{{\bf h}^{(2)}}^n\delta_{{\bf h}^{(2)}}\ldots s_{{\bf h}^{(\beta)}}^n\delta_{{\bf h}^{(\beta)}},  & & {\bf h}^{(2)}+\ldots+{\bf h}^{(\beta)}=
                                                                                                                   {\bf m}^{(1)}-{\bf m}^{(2)},\\
\hdots& &\hdots\hdots\hdots\hdots& &\hdots\hdots\hdots\hdots\hdots\hdots\hdots\hdots\\
\Delta_r&=&s_{{\bf j}^{(1)}}^n\delta_{{\bf j}^{(1)}}\ldots s_{{\bf j}^{(\tau)}}^n\delta_{{\bf j}^{(\tau)}}, & & {\bf j}^{(1)}+\ldots+{\bf j}^{(\tau)}={\bf m}^{(r)}, \quad \tau\geqslant2,
\end{array}
$$
the modules of all multi-indices in each $\delta$ being $\leqslant|{\bf m}|/2$.

Using \eqref{lemma3} for the $\tau$ single factors of $\Delta_r$ and applying \eqref{estimate} for the estimation of $\Delta_p$, $p=0,\,1,\ldots,\,r-1$, we find the inequalities
$$
\Delta_r\leqslant \left(|{\bf j}^{(1)}|\ldots|{\bf j}^{(\tau)}|\right)^{-2\nu n}N_2^{|{\bf m}^{(r)}|-\tau}Q^{|{\bf m}^{(r)}|},
$$
$$
\Delta_p\leqslant \left(|{\bf m}^{(p)}|-|{\bf m}^{(p+1)}|\right)^{-2\nu n}N_2^{|{\bf m}^{(p)}|-|{\bf m}^{(p+1)}|-1}Q^{|{\bf m}^{(p)}|-|{\bf m}^{(p+1)}|}, \qquad p=0,\,1,\ldots,\,r-1,
$$
whence
$$
\Delta_0\ldots\Delta_r\leqslant N_2^{|{\bf m}|-r-\tau}Q^{|{\bf m}|}\Bigl(\prod\limits_{i=1}^\tau|{\bf j}^{(i)}|\,\prod\limits_{p=0}^{r-1}
(|{\bf m}^{(p)}|-|{\bf m}^{(p+1)}|)\Bigr)^{-2\nu n}.
$$
By Corollary \ref{cor1},
$$
s_{{\bf m}^{(0)}}^n\epsilon_{{\bf m}^{(0)}}\ldots s_{{\bf m}^{(r)}}^n\epsilon_{{\bf m}^{(r)}}\leqslant N_1^{n(r+1)}|{\bf m}|^{\nu n}\,
\prod\limits_{p=0}^{r-1}(|{\bf m}^{(p)}|-|{\bf m}^{(p+1)}|)^{\nu n},
$$
and consequently from (\ref{decomp}) we get
$$
s_{\bf m}^n\delta_{\bf m}\leqslant N_1^{n(r+1)}N_2^{|{\bf m}|-r-\tau}Q^{|{\bf m}|} \Bigl(|{\bf m}|^{-1}\,\prod\limits_{i=1}^\tau|{\bf j}^{(i)}|^2\, \prod\limits_{p=0}^{r-1}(|{\bf m}^{(p)}|-|{\bf m}^{(p+1)}|)\Bigr)^{-\nu n}=
$$
$$
=N_1^{n(r+1)}N_2^{|{\bf m}|-r-\tau}Q^{|{\bf m}|} \Bigl(|{\bf m}|^{-1}\,\prod\limits_{i=1}^\tau y_i^2\,\prod\limits_{p=0}^{r-1}x_p\Bigr)^{-\nu n},
$$
where $x_p=|{\bf m}^{(p)}|-|{\bf m}^{(p+1)}|$, $y_i=|{\bf j}^{(i)}|$. Since $\sum\limits_{p=0}^{r-1}x_p+\sum\limits_{i=1}^\tau y_i=|{\bf m}|$, $\sum\limits_{i=1}^\tau y_i>|{\bf m}|/2$ and each $y_i\leqslant|{\bf m}|/2$, by Siegel's Lemma 1 \cite{Siegel} one has
$$
\prod\limits_{i=1}^\tau y_i^2\,\prod\limits_{p=0}^{r-1}x_p\geqslant \frac{|{\bf m}|^3}{8^{r+\tau-1}}.
$$
Therefore we come to the required estimate
$$
s_{\bf m}^n\delta_{\bf m}\leqslant N_1^{n(r+1)}N_2^{|{\bf m}|-r-\tau}Q^{|{\bf m}|}|{\bf m}|^{-2\nu n}\;8^{(r+\tau-1)\nu n}\leqslant |{\bf m}|^{-2\nu n}
N_2^{|{\bf m}|-1}Q^{|{\bf m}|}\frac{N_1^{n(r+\tau-1)}}{N_2^{r+\tau-1}}\,8^{(r+\tau-1)\nu n}=
$$
$$
=|{\bf m}|^{-2\nu n}N_2^{|{\bf m}|-1}Q^{|{\bf m}|}.
$$

\subsection{Finishing the proof}

By Lemmas \ref{lemmamaj}, \ref{l3} we get that the $s$-variate power series (\ref{formsol}) has a nonempty polydisc of convergence, say, ${\rm D}=\{|z_1|\leqslant\rho,\ldots,|z_s|\leqslant\rho\}$. Since all ${\rm Re}\,\alpha_1,\ldots,{\rm Re}\,\alpha_s>0$, taking any sector $S\subset\mathbb C$ with the vertex at the origin and of the opening less than $2\pi$, of such small radius that if $z\in S$ then $(z_1,\ldots,z_s)=(z^{\alpha_1},\ldots,z^{\alpha_s})\in\rm D$,
one concludes the proof of Theorem \ref{th1} by the equality
$$
\psi=\sum_{|{\bf m}|>0}c_{\bf m}\,z^{m_1\alpha_1+\ldots+m_s\alpha_s}=\sum_{|{\bf m}|>0}c_{\bf m}\,z_1^{m_1}\ldots z_s^{m_s}|_{z_i=z^{\alpha_i}}.
$$

\section{Theorem 1 bis and examples}

We note that the situation where all $q_i=q^{\alpha_i}$ lie {\it strictly inside} or {\it strictly outside} the unit circle (or equivalently, where all the $\alpha_i$'s lie {\it strictly above} or {\it strictly under} the line $\cal L$ passing through $0\in\mathbb C$ and having the slope $\ln |q|/\arg q$), which was considered in our previous paper \cite{GGL}, may be considered as a particular case of Theorem \ref{th1} (the basic estimates of Lemma \ref{l0} are fulfilled). At the same time taking into account Remark 2 we see that some other particular placements of the $q_i$'s with respect to the unit circle (of the $\alpha_i$'s with respect to the line $\cal L$) allow one to weak assumptions of Theorem \ref{th1}. Therefore we formulate a separate statement which follows from Theorem \ref{th1} and distinguishes all these particular cases of the placement of the $\alpha_i$'s on the plane.
\medskip

\noindent{\bf Theorem 1 bis} {\it The statement of Theorem \ref{th1} holds in the following particular cases:
\smallskip

{\rm a)} $L(0)\ne0$ and all the $\alpha_i$'s lie {\rm strictly above} the line $\cal L$;

{\rm b)} $\deg L=n$ and all the $\alpha_i$'s lie {\rm strictly under} the line $\cal L$;

{\rm c)} all the $\alpha_i$'s lie {\rm on} the line $\cal L$ and the condition \eqref{smalldiv+} is fulfilled for those roots $\xi=a$ of the polynomial $(\xi-q^{\lambda_{j_0}})L(\xi)$ that lie on the circle $\{|\xi|=|q^{\lambda_{j_0}}|\}$;

{\rm d)} $L(0)\ne0$, all the $\alpha_i$'s lie {\rm above} or {\rm on} the line $\cal L$, and the condition \eqref{smalldiv+} is fulfilled for those roots $\xi=a$ of the polynomial $(\xi-q^{\lambda_{j_0}})L(\xi)$ that lie inside the closed disk $\{|\xi|\leqslant|q^{\lambda_{j_0}}|\}$;

{\rm e)} $\deg L=n$, all the $\alpha_i$'s lie {\rm under} or {\rm on} the line $\cal L$, and the condition \eqref{smalldiv+} is fulfilled for those roots $\xi=a$ of the polynomial $(\xi-q^{\lambda_{j_0}})L(\xi)$ that lie outside the open disk $\{|\xi|<|q^{\lambda_{j_0}}|\}$.}
\medskip

Note that the small divisors phenomenon for classical power series solutions of (\ref{e1}) arising in the case of $q=e^{2\pi{\rm i}\omega}$, $\omega\in{\mathbb R} \setminus\mathbb Q$, and studied in \cite{b1}, \cite{dv}, is contained in the case c) of Theorem 1 bis: the line $\cal L$ coincides with the $Ox$ axis, $\lambda_{j_0}=0$, the set of power exponents is generated by the unique $\alpha_1=1\in\cal L$ and the condition (\ref{smalldiv+}) is reduced to (\ref{e2equiv}) in this case. 

Further we consider several examples.
\medskip

{\bf Example 1.} We present an example that expectingly illustrates an importance of the assumptions $L(0)\ne0$, $\deg L=n$ of Theorem \ref{th1}. More precisely, we consider a situation where the $\alpha_i$'s are placed according to the case d) of Theorem 1 bis and satisfy (\ref{smalldiv+}) but $L(0)=0$ and a corresponding generalized power series solution diverges.

An equation
$$
\sigma^2 y-q^{\rm i}\,\sigma y+q^{2\rm i}\,z(1+y^2)=0, \qquad q=e^{-\frac{\sqrt2}2\pi(1-\rm i)},
$$
possesses a generalized formal power series solution
$$
\varphi=c_{0,0}\,z^{\rm i}+\sum\limits_{m_1+m_2>0}c_{m_1,m_2}\,z^{{\rm i}+m_1(1+{\rm i})+m_2(1-{\rm i})},
$$
where the complex coefficient $c_{0,0}\neq 0$ is arbitrary and the other complex coefficients $c_{m_1,m_2}$ are determined uniquely by $c_{0,0}$. Indeed, substituting $y=c_0\,z^{\lambda}+\ldots$ into the equation one obtains
$$
c_0(q^{2\lambda}-q^{\lambda+\rm i})z^{\lambda}+\ldots=-q^{2\rm i}\,z+\ldots,
$$
and since $q^2-q^{1+\rm i}\ne-q^{2\rm i}$, this implies $q^{2\lambda}-q^{\lambda+\rm i}=0$, that is, $\lambda=\rm i$. Further making the change of the unknown
$y=c_{0,0}\,z^{\rm i}+z^{\rm i}\,u$ one comes to an equation
\begin{equation}\label{eqex1}
\sigma^2 u-\sigma u=-z^{1-\rm i}-z^{1+\rm i}(c_{0,0}+u)^2,
\end{equation}
whence the form of the generalized power series solution $\varphi$ follows.

Note that
\begin{eqnarray*}
\alpha_1=1+\rm i & \mbox{and} & q_1=q^{\alpha_1}=e^{-\sqrt2\pi}<1; \\
\alpha_2=1-\rm i & \mbox{and} & q_2=q^{\alpha_2}=e^{\sqrt2\pi\rm i}, \; |q_2|=1.
\end{eqnarray*}
Thus concerning the unique nonzero root $a=q^{\rm i}$ of the polynomial $L(\xi)=\xi(\xi-q^{\rm i})$, the diophantine condition (\ref{smalldiv+}) is fulfilled for it. Indeed, $\lambda_{j_0}=\rm i$ and the left hand side of (\ref{smalldiv+}) becomes
$$
\Bigl|\bigl(m_1(1+{\rm i})+m_2(1-{\rm i})\bigr)\frac{\sqrt2}2\pi(1-{\rm i})+2\pi m{\rm i}\Bigr|=\pi\bigl|m_1\sqrt2+2{\rm i}(m-m_2\sqrt2/2)\bigr|,
$$
which is not less than
$$
{\rm const}\cdot|m-m_2\sqrt2/2|>cm_2^{-\nu}\geqslant c(m_1+m_2)^{-\nu}
$$
for some positive constants $c,\nu$, in view of the algebraicity of the irrational number $\sqrt2/2$.

Theorem 1 cannot be applied since $L(0)=0$. Let us show that $\varphi$ has a zero radius of convergence. Consider a subseries
$$
\psi_+=\sum_{m=1}^\infty c_{m,0}\,z^{m(1+{\rm i})}
$$
of the generalized power series $\psi=\sum_{m_1+m_2>0}c_{m_1,m_2}\,z^{m_1(1+{\rm i})+m_2(1-{\rm i})}$. Since the latter satisfies (\ref{eqex1}), this $\psi_+$ satisfies an equation
$$
\sigma^2 u-\sigma u=-z^{1+\rm i}(c_{0,0}+u)^2,
$$
whence follows an equality
$$
\sum_{m=1}^\infty q^{m(1+{\rm i})}(q^{m(1+{\rm i})}-1)c_{m,0}\,z^{m(1+{\rm i})}=-c_{0,0}^2\,z^{1+\rm i}-\sum_{m=2}^\infty2c_{0,0}c_{m-1,0}\,z^{m(1+{\rm i})}-
z^{1+\rm i}\Bigl(\sum_{m=1}^\infty c_{m,0}\,z^{m(1+{\rm i})}\Bigr)^2.
$$
Thus we have the following recurrence relations for the coefficients $c_{m,0}$:
$$
c_{1,0}=\frac{c_{0,0}^2}{e^{-\sqrt2\pi}(1-e^{-\sqrt2\pi})}, \qquad c_{m,0}=\frac{2c_{0,0}c_{m-1,0}+\sum_{j=1}^{m-2}c_{j,0}c_{m-1-j,0}}
{e^{-\sqrt2\pi m}(1-e^{-\sqrt2\pi m})}, \quad m\geqslant2.
$$
If $c_{0,0}$ is a positive real number, the same all the $c_{m,0}$'s are. Taking $c_{0,0}\geqslant1/2$ we obtain
$$
c_{m,0}\geqslant e^{\sqrt2\pi m}c_{m-1,0}\geqslant\ldots\geqslant e^{\sqrt2\pi m(m+1)/2}c_{0,0},
$$
whence the divergence of $\psi_+$ and hence of $\psi$ and $\varphi$ follows.
\medskip

{\bf Example 2.} Consider a kind of a $q$-difference analogue of the Painlev\'e III equation with $a=b=0$, $c=d=1$:
$$
y\,\sigma^2 y-(\sigma y)^2-z^2y^4-z^2=0,
$$
where $q=e^{2{\rm i}\pi\omega}$, $\omega\in\mathbb{R}\setminus\mathbb{Q}$.
This equation possesses a two-parameter family of formal solutions:
$$
\varphi=\sum_{m_1,m_2\in{\mathbb Z}_+}c_{m_1,m_2}\,z^{r+m_1(2-2r)+m_2(2+2r)},
$$
where the complex coefficient $c_{0,0}\ne0$ is arbitrary, $-1<{\rm Re}\,r<1$, the other complex coefficients $c_{m_1,m_2}$ are uniquely determined by $c_{0,0}$ and $r$. Since the numbers $q^{2\pm2r}$ lie on the opposite sides of the unit circle (if ${\rm Im}\,r\ne0$) or on the unit circle (if ${\rm Im}\,r=0$), we cannot apply results of our previous paper \cite{GGL}. However, the polynomial $L(\xi)=c_{0,0}(\xi-q^r)^2$ does not vanish at $0$, therefore taking $r$, $\omega$ in such a way that the condition of Theorem \ref{th1} holds,
$$
\left|(m_1(2-2r)\omega+m_2(2+2r)\omega-m\right|>c\,|m_1+m_2|^{-\nu}
$$
for some positive $c$ and $\nu$, we obtain the convergent $\varphi$. For example, it is sufficient for $\omega$ to be algebraic and for $r$ simply to have a non-zero imaginary part. 
\medskip

{\bf Example 3.} Let $$L(\xi)=(\xi-2)(\xi-{\rm i}^{4\ell})(\xi-{\rm i}^{1+{\rm i}}),$$ where
$$\displaystyle \ell=\sum\limits_{k=1}^{\infty}10^{-k!}
$$ is the Liouville number, and let us consider a $q$-difference equation
$$  L(\sigma)\,y-z^7y^2-z^2=0, \qquad q={\rm i},
$$
which possesses a formal solution
$$
\varphi=z^{4\ell}+z^{1+{\rm i}}+\sum\limits_{m_1+m_2+m_3\geqslant 2}c_{m_1,m_2,m_3}\,z^{m_1+m_2\,4\ell+m_3(1+{\rm i})}, \quad c_{m_1,m_2,m_3}\in\mathbb{C}.
$$
Here $\deg L=3$, $L(0)\neq 0$, but the arithmetic condition of Theorem \ref{th1} is not fulfilled for the root $a={\rm i}^{4\ell}$ of the polynomial $L$. Indeed, since $\lambda_{j_0}$ belongs to the set of power exponents generated by $1, 4\ell, 1+\rm i$ in this case, the violence of (\ref{smalldiv+}) means that for every $\nu$, $c>0$ there exist integers $m$ and $m_1,m_2,m_3\geqslant 0$,  $m_1+m_2+m_3\geqslant 2$, such that the inequality
\begin{equation}\label{t1}
\left|\left(m_1+m_2\,4\ell+m_3(1+{\rm i})\right)\ln {\rm i}-\ln {\rm i}^{4\ell}-2\pi m{\rm i}\right|>c\,|m_1+m_2+m_3|^{-\nu}
\end{equation}
does not hold. Let us prove this. In the case $m_1=m_3=0,$ $m_2=p+1$, $p\in\mathbb{N},$ the inequality \eqref{t1} becomes
\begin{equation}\label{t2}
\left|(p+1)\,4\ell\;\ln {\rm i}-\ln{{\rm i}^{4\ell}}-2\pi m {\rm i}\right|> c\,(p+1)^{-\nu}.
\end{equation}
The latter can be rewritten in a simpler form
\begin{equation}\label{t3}
\left|p\,\ell-m\right|>c\,(p+1)^{-\nu}.
\end{equation}
But the inequality \eqref{t3} is violated for infinitely many pairs $(p,m)$ whatever $\nu$, $c>0$ we take. Indeed, if we take  $\tilde{p}=10^{N!}$, $\tilde{m}=\sum\limits_{k=1}^{N}10^{N!-k!}$, then
$$
\left|\tilde{p}\,\ell-\tilde{m}\right|=\sum\limits_{k=N+1}^\infty 10^{N!-k!}=\sum\limits_{j=0}^\infty\tilde{p}^{1-\prod\limits_{i=0}^j (N+1+i)}< \tilde{p}^{\,-N+1}<c\,(\tilde{p}+1)^{-\nu}
$$
for $N$ big enough.
Therefore,  in this example we cannot guarantee the convergence of the formal solution $\varphi$.

\end{document}